\documentclass{elsart}

\usepackage{epsfig}
\usepackage{amsfonts,amsmath}
\usepackage[latin1]{inputenc}
\newcounter{subfigure}

\begin{document}

\begin{frontmatter}

\title{Approximation of the domain of attraction of an asymptotically stable fixed point of a first order
analytical system of difference equations}
\author[math]{E. Kaslik}
\author[physics]{A.M. Balint}
\author[math]{S. Birauas}
\author[math]{St. Balint}
\address[math]{Department of Mathematics,\hfill\break West
University of Timisoara, Bd V. Parvan No. 4, 1900 Timisoara,
Romania}
\address[physics]{Department of Physics,\hfill\break West
University of Timisoara, Bd V. Parvan No. 4, 1900 Timisoara,
Romania}

\begin{abstract}
In this paper a first order analytical system of difference
equations is considered. For an asymptotically stable fixed point
$x^{0}$ of the system a gradual approximation of the domain of
attraction ($DA$) is presented in the case when the matrix of the
linearized system in $x^{0}$ is a contraction. This technique is
based on the gradual extension of the "embryo" of an analytic
function of several variables. The analytic function is a Lyapunov
function whose natural domain of analyticity is the $DA$ and which
satisfies an iterative functional equation. The equation permits
to establish an "embryo" of the Lyapunov function and a first
approximation of the $DA$. The "embryo" is used for the
determination of a new "embryo" and a new part of the $DA$. In
this way, computing new "embryos" and new domains, the $DA$ is
gradually approximated. Numerical examples are given for
polynomial systems.
\end{abstract}


\end{frontmatter}


\section{Introduction}

We consider the system of difference equations

\begin{equation}
\label{ec} x_{k+1}=g(x_{k})\qquad k=0,1,2...
\end{equation}
where $g:\Omega\rightarrow\Omega$ is an analytic function defined
on a domain $\Omega$ included in $\mathbb{R}^{n}$.

A point $x^{0}\in\Omega$ is a fixed point for the system
(\ref{ec}) if $x^{0}$ satisfies

\begin{equation}\label{fix}
    x^{0}=g(x^{0})
\end{equation}

The fixed point $x^{0}$ of (\ref{ec}) is "stable" provided that
given any ball
$B(x^{0},\varepsilon)=\{x\in\Omega/\|x-x^{0}\|<\varepsilon\}$,
there is a ball
$B(x^{0},\delta)=\{x\in\Omega/\|x-x^{0}\|<\delta\}$ such that if
$x\in B(x^{0},\delta)$ then $g^{k}(x)\in B(x^{0},\varepsilon)$,
for $k=0,1,2,...$ [1].

If in addition there is a ball $B(x^{0},r)$ such that
$g^{k}(x)\rightarrow x^{0}$ as $k\rightarrow\infty$ for all $x\in
B(x^{0},r)$ then the fixed point $x^{0}$ is "asymptotically
stable".[1].

The domain of attraction $DA(x^{0})$ of the asymptotically stable
fixed point $x^{0}$ is the set of initial states $x\in \Omega$
from which the system converges to the fixed point itself i.e.

\begin{equation}\label{da}
    DA(x^{0})=\{x\in\Omega | g^{k}(x)\stackrel{k\rightarrow\infty}{\longrightarrow}x^{0}\}
\end{equation}

It is known that $x^{0}$ is a fixed point for system (\ref{ec}) if
and only if $0\in\mathbb{R}^{n}$ is a fixed point for the system

\begin{equation}
\label{ec0} y_{k+1}=f(y_{k}) \qquad k=0,1,2...
\end{equation}
where $f:\Omega-x^{0}\rightarrow\Omega-x^{0}$ is the analytic
function defined by

\begin{equation}\label{f}
    f(y)=g(y+x^{0})-x^{0} \qquad \textrm{for } y\in\Omega-x^{0}
\end{equation}

The fixed point $x^{0}$ of (\ref{ec}) is asymptotically stable if
and only if the fixed point $0\in\mathbb{R}^{n}$ of the system
(\ref{ec0}) is asymptotically stable.

The domain of attraction of $x^{0}$, $DA(x^{0})$ is related to the
domain of attraction of $0$, $DA(0)$ by the equation

\begin{equation}\label{relda}
    DA(x^{0})=DA(0)+x^{0}
\end{equation}

For the above reason in the followings instead of the system
(\ref{ec}) we will consider the system (\ref{ec0}).

Theoretical research shows that the $DA(0)$ and its boundary are
complicated sets [2],[3],[4],[5],[6]. In most cases, they do not
admit an explicit elementary representation. For this reason,
different procedures are used for the approximation of the $DA(0)$
with domain having a simpler shape. For example, in the case of
the theorem 4.20 pg 170 [1] the domain which approximates the
$DA(0)$ is defined by a Lyapunov function $V$ built with the
matrix $\partial_{0}f$ of the linearized system in $0$. In this
paper, we present a technique for the construction of a Lyapnov
function $V$ in the case when the matrix $\partial_{0}f$ is a
contraction, i.e. $\|\partial_{0}f\|<1$. The Lyapunov function $V$
is built using the whole nonlinear system, not only the matrix
$\partial_{0}f$. $V$ is defined on the whole $DA(0)$, and more,
the $DA(0)$ is the natural domain of analyticity of $V$. The
formula which defines the Lyapunov function $V$ is used for
determining an "embryo" of $V$ and a first approximation of
$DA(0)$. The "embryo" is used for the determination of a new
"embryo" and a new part of $DA(0)$. In this way, computing new
"embryos" and new domains the $DA(0)$ is gradually approximated.

\section{Theoretical results}

Let be $f:\Omega\rightarrow\Omega$ an analytic function defined on
a domain $\Omega\subset\mathbb{R}^{n}$ containing the origin
$0\in\mathbb{R}^{n}$.

\begin{thm}
If the funtion $f$ satisfies the following conditions:
\begin{equation}
  f(0) = 0
\end{equation}
\begin{equation}
  \|\partial_{0}f\| < 1
\end{equation}
then $0$ is an asymptotically stable fixed point. $DA(0)$ is an
open subset of $\Omega$ and coincides with the natural domain of
analyticity of the unique solution $V$ of the iterative first
order functional equation
\begin{equation}
\label{ecV}
\begin{array}{ll}
\left\{\begin{array}{l}
V(f(x))-V(x)=-\|x\|^{2}\\
V(0)=0
\end{array}\right.
\end{array}
\end{equation}
The function $V$ is positive on $DA(0)$ and
$V(x)\stackrel{x\rightarrow x^{0}}{\longrightarrow}+\infty$, for
any $x^{0}\in FrDA(0)$ ($FrDA(0)$ denotes the boundary of
$DA(0)$).
\end{thm}

\begin{pf} Let be $\alpha$ such that $\|\partial_{0} f\|<\alpha<1$. By
the continuity of $x\mapsto\|\partial_{x} f\|$ there exists a
$\delta>0$ such that $\|\partial_{x} f\|\leq\alpha$ for
$\|x\|\leq\delta$. The mean value theorem gives

\begin{equation}
    \label{t1}
    \|f(x')-f(x'')\|\leq\alpha\|x'-x''\|
\end{equation}
for any $x'$ and $x''$ in the ball $B(0,\delta)$. Therefore
\begin{equation}
    \label{t2}
    \|f^{k}(x)\|\leq\alpha^{k}\|x\|
\end{equation}
for any $x$ in the ball $B(0,\delta)$ and $k=0,1,2,...$. For a
ball $B(0,\varepsilon)$ we take $\delta'=\min(\varepsilon,\delta)$
and the ball $B(0,\delta')$. We have $f^{k}(x)\in
B(0,\varepsilon)$ for any $x\in B(0,\delta')$ and $k=0,1,2...$,
which means that $0$ is stable.

From (\ref{t2}) we obtain
$f^{k}(x)\stackrel{k\rightarrow\infty}{\longrightarrow}0$ for any
$x\in B(0,\delta)$. We can conclude now that $0$ is asymptotically
stable.

In order to show that $DA(0)$ is an open subset of $\Omega$ we
consider $x'$ from $DA(0)$ and $k_{0}>0$ such that
$\|f^{k_{0}}(x')\|<\frac{\delta}{3}$. Because $f^{k_{0}}$ is a
continuous function, there exists a ball $B(x',\delta'')$ such
that $\|f^{k_{0}}(x)-f^{k_{0}}(x')\|<\frac{\delta}{3}$, for any
$x\in B(x',\delta'')$. Therefore,
$\|f^{k_{0}}(x)\|\leq\frac{2\delta}{3}<\delta$, for any $x\in
B(x',\delta'')$. It follows that $x\in DA(0)$ and therefore,
$B(x',\delta'')\subset DA(0)$. This proves that $DA(0)$ is an open
subset of $\Omega$.

Now we consider $x\in DA(0)$ and the sequence
$\{f^{k}(x)\}_{k\in\mathbb{N}}$. There exists $k_{x}$ such that
$f^{k}(x)\in B(0,\delta)$, for any $k\geq k_{x}$. Therefore,
$\|f^{k_{x}+k}(x)\|\leq\alpha^{k}\|f^{k_{x}}(x)\|$, for
$k=0,1,2...$. It follows that the series
$\sum\limits_{k=0}^{\infty}\|f^{k}(x)\|^{2}$ is convergent for any
$x\in DA(0)$.

Let be $V=V(x)$ the function defined by

\begin{equation}\label{t3}
    V(x)=\sum\limits_{k=0}^{\infty}\|f^{k}(x)\|^{2} \qquad \textrm{ for }
    x\in DA(0)
\end{equation}

The above function defined on $DA(0)$ is analytical, positive and
satisfies (\ref{ecV}). In order to show that the function $V$
defined by (\ref{t3}) is the unique function which satisfies
(\ref{ecV}) we consider $V'=V'(x)$ satisfying (\ref{ecV}) and we
denote by $V''$ the difference $V''=V-V'$. It is easy to see that
$V''(f(x))-V''(x)=0$, for any $x\in DA(0)$. Therefore, we have
$V''(x)=V''(f^{k}(x))$ for any $x\in DA(0)$ and any $k=0,1,2...$.
It follows that
$V''(x)=\lim\limits_{k\rightarrow\infty}V''(f^{k}(x))=0$ for any
$x\in DA(0)$. In other words, $V(x)=V'(x)$, for any $x\in DA(0)$,
so $V$ defined by (\ref{t3}) is the unique function which
satisfies (\ref{ecV}).

In order to show that $V(x)\stackrel{x\rightarrow
x^{0}}{\longrightarrow}\infty$ for any $x^{0}\in  FrDA(0)$ we
consider $x^{0}\in  FrDA(0)$ and $r>0$ such that
$\|f^{k}(x^{0})\|>r$, for any $k=0,1,2...$. For an arbitrary
positive number $N>0$ we consider the first natural number $k_{1}$
which satisfies $k_{1}\geq\frac{2N}{r^{2}}+1$. Let be $r_{1}>0$
such that $\|f^{k}(x)\|\geq\frac{r}{\sqrt{2}}$ for any
$k=1,2,..,k_{1}$ and $x\in B(x^{0},r_{1})$. For any $x\in
B(x^{0},r_{1})\cap DA(0)$ we have
$\sum\limits_{k=0}^{k_{1}}\|f^{k}(x)\|^{2}>N$. Therefore,
$V(x)\stackrel{x\rightarrow x^{0}}{\longrightarrow}\infty$.\qed
\end{pf}

\begin{rem}
Newton's method for solving systems of $n$ nonlinear equations
with $n$ unknowns always leads to a system of difference equations
which satisfies the conditions of Theorem 1.
\end{rem}

We consider the expansion of $V$ in $0$:
\begin{equation}\label{serieV}
    V(x_{1},x_{2},...,x_{n})=\sum\limits_{m=2}^{\infty}\sum\limits_{|j|=m}B_{j_{1}j_{2}...j_{n}}x_{1}^{j_{1}}x_{2}^{j_{2}}...x_{n}^{j_{n}}
\end{equation}

The domain of convergence of the series (\ref{serieV})is the set
of those $x\in\Omega$ which have the property that the series
(\ref{serieV}) is absolutely convergent in a neighborhood of $x$
[7]. We denote by $D^{0}$ this domain. Actually, the domain of
convergence $D^{0}$ coincides with the interior of the set of all
points $x^{0}$ in which the series (\ref{serieV}) is absolutely
convergent [7].

\begin{thm}
(Cauchy-Hadamard see [7]) A point $x$ belongs to $D^{0}$ if and
only if
\begin{equation} \label{ineg}
\overline{\lim_{m}}\sqrt[m]{\sum_{|j|=m}|B_{j_{1}j_{2}...j_{n}}x_{1}^{j_{1}}x_{2}^{j_{2}}...x_{n}^{j_{n}}|}<1
\end{equation}
\end{thm}

By the previous theorem we find a part $D^{0}$ of $DA(0)$.

In practice, we can compute the coefficients
$B_{j_{1}j_{2}...j_{n}}$ up to a finite degree $p$. This degree
$p$ has to be big enough for assuring that the domain $D^{0}_{p}$
given by

\begin{equation}
D^{0}_{p}=\{x\in\Omega/\sqrt[p]{\sum\limits_{|j|=p}|B_{j_{1}j_{2}...j_{n}}
 x_{1}^{j_{1}}x_{2}^{j_{2}}...x_{n}^{j_{n}}|}<1\}
 \nonumber
\end{equation}
\noindent approximates the region of convergence $D^{0}$ of the
series of $V$ and the "embryo"
\begin{equation}
 V^{0}_{p}(x_{1},x_{2},...,x_{n})=\sum\limits_{m=2}^{p}\sum\limits_{|j|=m}B_{j_{1}j_{2}...j_{n}}
 x_{1}^{j_{1}}x_{2}^{j_{2}}...x_{n}^{j_{n}}
\end{equation}
approximates $V$ with accuracy.

The first estimate of the $DA$ will be $D^{0}_{p}$.

In order to extend the first estimate $D^{0}_{p}$, we expand $V$
in a point $x^{0}$ close to $ FrD^{0}_{p}$ in which
$|V^{0}_{p}(x^{0})|$ is still small. That is because, according to
Theorem 1, the points $x$ close to $ FrD^{0}_{p}$ for which
$|V^{0}_{p}(x)|$ is extremely high are close to $ FrDA(0)$.

To find the expansion of $V$ in $x^{0}$ close to $ FrD^{0}_{p}$,
we will compute the expansion in $x^{0}$ of the "embryo"
$V^{0}_{p}$ of $V$. We obtain:

\begin{eqnarray}
V^{1}_{p}(x_{1},x_{2},...,x_{n})&=&
\sum\limits_{m=0}^{p}\sum\limits_{|j|=m}\frac{\partial^{m}V^{0}_{p}}{\partial
x_{1}^{j_{1}}\partial x_{2}^{j_{2}}...\partial
x_{n}^{j_{n}}}(x^{0})
 (x_{1}-x_{1}^{0})^{j_{1}}...(x_{n}-x_{n}^{0})^{j_{n}}
 /m!=\nonumber\\
 &=&
 \sum\limits_{m=0}^{p}\sum\limits_{|j|=m}B^{1}_{j_{1}j_{2}...j_{n}}
 (x_{1}-x_{1}^{0})^{j_{1}}(x_{2}-x_{2}^{0})^{j_{2}}...(x_{n}-x_{n}^{0})^{j_{n}}
\end{eqnarray}

\noindent We consider the set

\begin{equation}
D_{p}^{1}=\{x\in\Omega/\sqrt[p]{\sum\limits_{|j|=p}|B_{j_{1}j_{2}...j_{n}}
 (x_{1}-x_{1}^{0})^{j_{1}}(x_{2}-x_{2}^{0})^{j_{2}}...(x_{n}-x_{n}^{0})^{j_{n}}|}<1\}
 \nonumber
\end{equation}

\noindent which provides a new part $D^{1}_{p}$ of $DA(0)$.

So, $D_{p}^{0}\cup D_{p}^{1}$ gives a larger estimate of $DA(0)$.
We can continue this procedure for a few steps, till the values
$|V_{p}^{k}|$ become extremely large and we obtain the estimate
$D_{p}^{0}\cup D_{p}^{1}\cup ... \cup D_{p}^{k}$ of $DA(0)$.

\section{Numerical results}

\subsection{Systems with known domains of attraction}

In this subsection, we will present some examples of systems of
one, two or three difference equations, for which we can compute
easily the $DA$. We will apply our technique to these examples,
and we will show how the real domains of attraction are gradually
approximated. These examples are meant to validate our procedure.

The computations were made using a program written in Mathematica
4, Wolfram Research on an Intel Pentium III PC (2Ghz, 512MB of
RAM). The data for the estimations (the degree up to which the
approximation is made, the necessary timing for the estimations)
are displayed in Table~1. In our figures, the thick black line
represents the true boundary of the domain of attraction, the dark
grey set denotes the first estimate $D_{p}^{0}$ of $DA$ and the
further estimates $D_{p}^{k}$ of $DA$ with $k\geq 1$ are colored
in light grey.

\newpage

\subsubsection*{Example 1}

We consider the following one dimensional difference equation:

\begin{equation}\label{ex1}
    x_{n+1}=\frac{1}{2}x_{n}-x_{n}^{2}+2x_{n}^{3}-4x_{n}^{4}
\end{equation}

This difference equation has two fixed points: $x=0$ and
$x=-0.271845$. The fixed point $x=0$ is asymptotically stable,
while the other fixed point is unstable. It can be proved (by the
staircase method) that $DA(0)$ is the interval
$(-0.271845,0.653564)$.

We applied the procedure described above to obtain an estimation
of the $DA(0)$.
\begin{itemize}
    \item After the first step, we obtained
    $D^{0}=(-0.27184,0.27184)$.
    \item After the second step applied in $x^{1}=0.2718$, we
    obtained $D^{1}=(0.01345,0.53015)$.
    \item After the third step applied in $x^{2}=0.5$, we
    obtained $D^{2}=(0.38378,0.61622)$.
    \item After the forth step applied in $x^{3}=0.61$, we
    obtained $D^{3}=(0.59785,0.622175)$.
\end{itemize}

Therefore, the estimate of $DA(0)$ obtained after four steps is
$(-0.27184,0.622175)$.

\subsubsection*{Example 2}

In [1], the following system of difference equations is
considered:

\begin{equation}
\label{ex2}
\begin{array}{ll}
\left\{\begin{array}{l}
x_{n+1}=-y_{n}(1-x_{n}^{2}-y_{n}^{2})\\
y_{n+1}=-x_{n}(1-x_{n}^{2}-y_{n}^{2})
\end{array}\right.
\end{array}
\end{equation}

\noindent It is easy to see that $(x,y)=(0,0)$ is an
asymptotically stable fixed point for the system (\ref{ex2}), and
the matrix of the linearized system in $(0,0)$ satisfies the
conditions from Theorem 1. There are two more fixed points of
system (\ref{ex2}) (which are unstable) which are represented by
the gray points in Fig 1.1-1.2. The $DA(0)$ is the set:
\begin{displaymath}
DA(0)=\{(x,y)\in\mathbb{R}^{2}/x^{2}+y^{2}<3\}
\end{displaymath}

\noindent Using the technique described above, after 2 steps of
gradual extension of the "embryo", we find the estimate of $DA(0)$
presented in Figure 1.1. If we apply the second step for several
points close to $ FrD^{0}$, we obtain the estimate of $DA(0)$
presented in Figure 1.2.

\renewcommand{\thefigure}{\arabic{figure}.\arabic{subfigure}}
\setcounter{subfigure}{1}

\begin{figure}[htbp]
\begin{minipage}[t]{0.5\linewidth}
\centering
\includegraphics*[bb=0cm 0cm 11cm 11cm,
width=5cm]{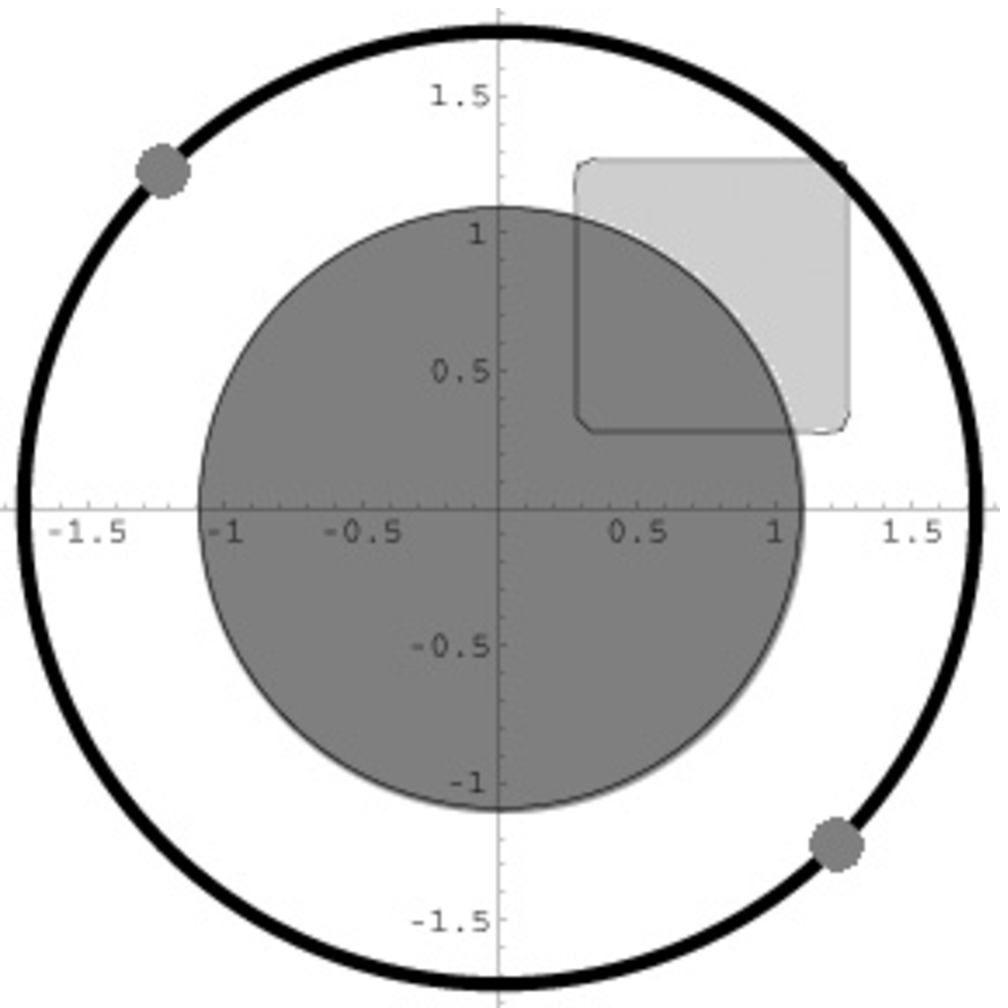} \caption{The estimate of $DA(0)$ after 2
steps for system (\ref{ex2})}
\end{minipage}
\addtocounter{figure}{-1} \addtocounter{subfigure}{1}
\begin{minipage}[t]{0.5\linewidth}
\centering
\includegraphics*[bb=0cm 0cm 11cm
11cm, width=5cm]{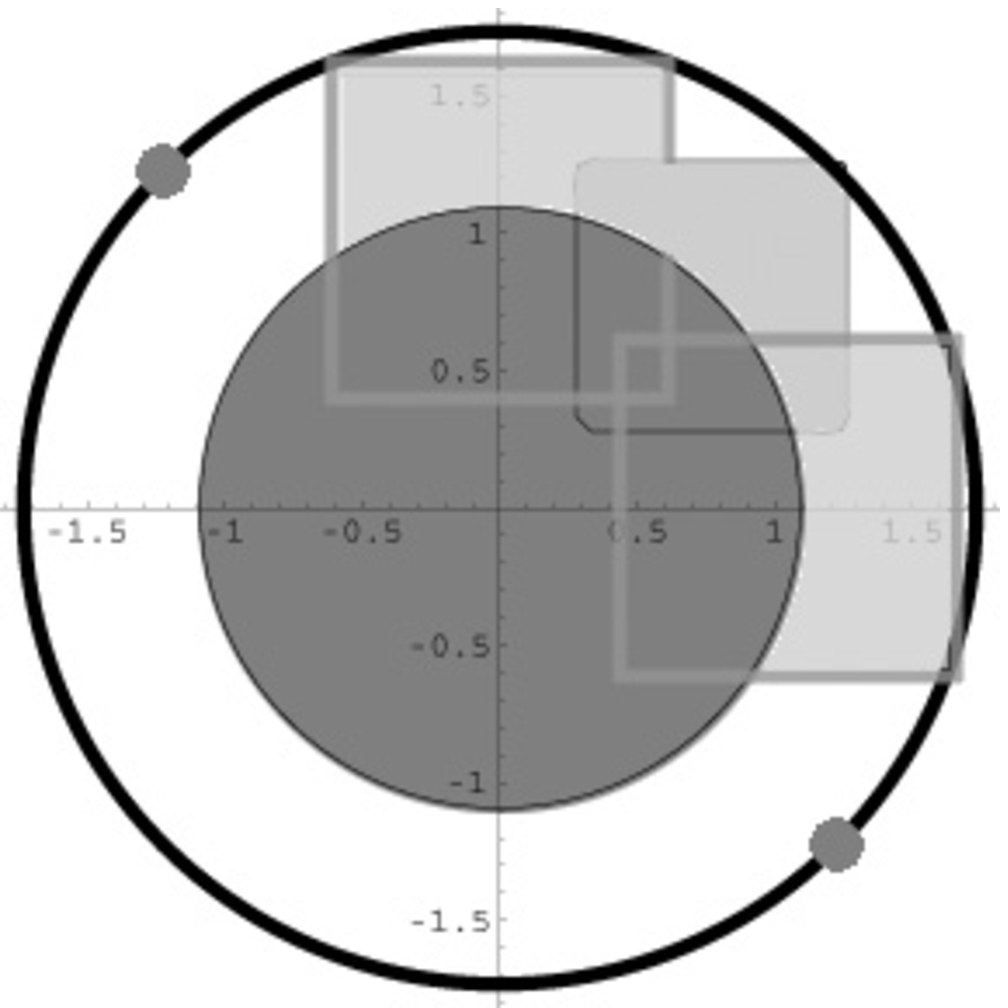} \caption{The estimate of the $DA(0)$
after 2 steps, in several points close to $ FrD^{0}$}
\end{minipage}
\end{figure}

\subsubsection*{Example 3}

The following decoupled system of difference equations is
considered:

\begin{equation}
\label{ex3}
\begin{array}{ll}
\left\{\begin{array}{l}
x_{n+1}=4x_{n}^{3}\\
y_{n+1}=9y_{n}^{3}
\end{array}\right.
\end{array}
\end{equation}

\noindent It is clear that $(0,0)$ is an asymptotically stable
fixed point for the system (\ref{ex3}) and
\begin{displaymath}
DA(0)=(-\frac{1}{2},\frac{1}{2})\times(-\frac{1}{3},\frac{1}{3})
\end{displaymath}

\noindent Applying the above procedure, after only one step, we
obtain the whole $DA(0)$.

\renewcommand{\thefigure}{\arabic{figure}}
\begin{figure}[htbp]
\centering
\includegraphics*[bb=0cm 0cm 11cm
11cm, width=5cm]{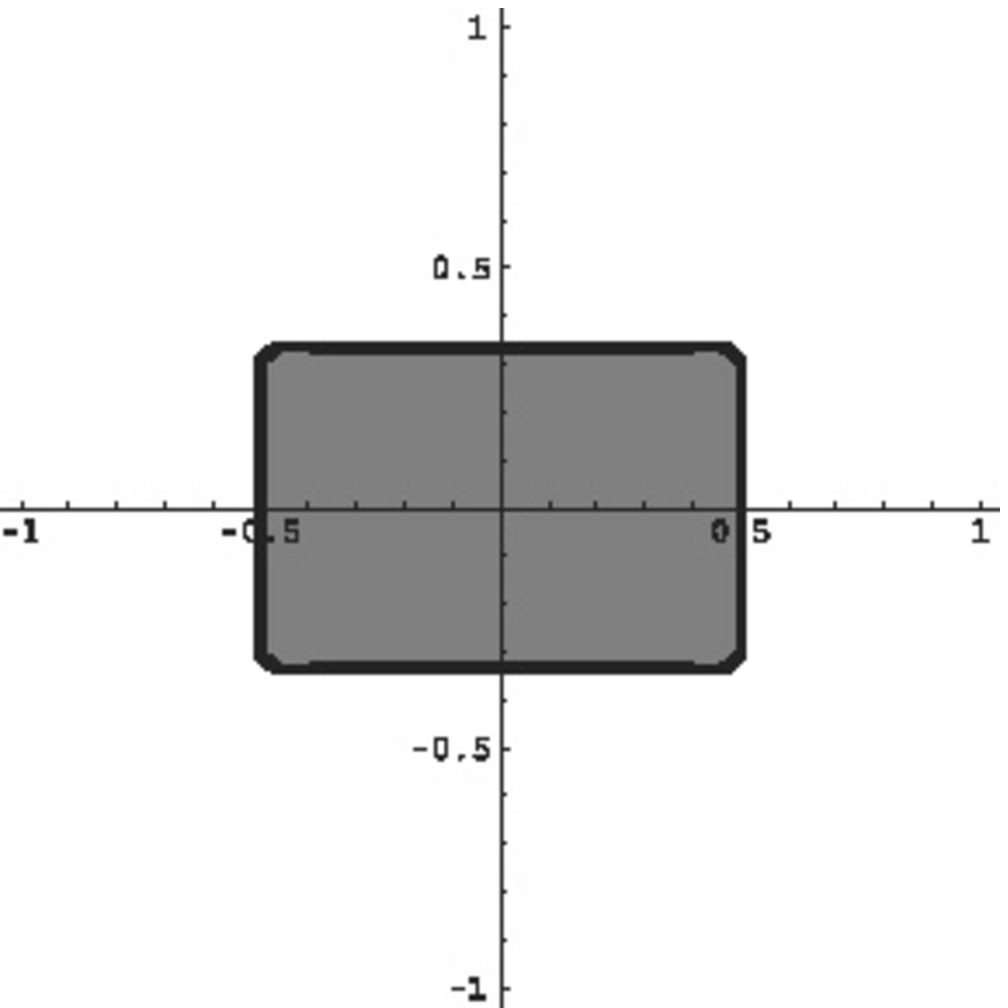} \caption{The estimate of $DA(0)$ after
1 steps for system (\ref{ex3})}
\end{figure}

\newpage

\subsection{Systems for which we don't know the domain of attraction}

In this subsection, some systems of difference equations are
presented for which we don't know the $DA(0)$. For these examples,
we will apply the procedure presented above.

\subsubsection*{Example 4}
We consider the difference equation:
\begin{equation}\label{ex4}
    x_{n+1}=-x_{n}^{2}-2x_{n}^{3}-4x_{n}^{4}-8x_{n}^{5}
\end{equation}

This equation has only one fixed point, namely $x=0$. We obtained
\begin{itemize}
    \item  $D^{0}=(-0.442585,0.442585)$, after the first step.
    \item  $D^{1}=(-0.673088,-0.212082)$, after the second step applied in $x^{1}=-0.44258$.
\end{itemize}

The estimate of the $DA(0)$ obtained after two steps is
$(-0.673088,0.442585)$.

\subsubsection*{Example 5}

The following system is considered:

\begin{equation}
\label{ex5}
\begin{array}{ll}
\left\{\begin{array}{l}
x_{n+1}=-\frac{1}{2}x_{n}+x_{n}y_{n}\\
y_{n+1}=-\frac{1}{2}y_{n}+x_{n}y_{n}
\end{array}\right.
\end{array}
\end{equation}

This system of difference equations has besides the asymptotically
stable fixed point $(0,0)$, an unstable fixed point represented in
the following figure by a black point. After two steps, we can
obtain the following estimate of the $DA(0)$:

\begin{figure}[htbp]
\centering
\includegraphics*[bb=0cm 0cm 11cm
11cm, width=5cm]{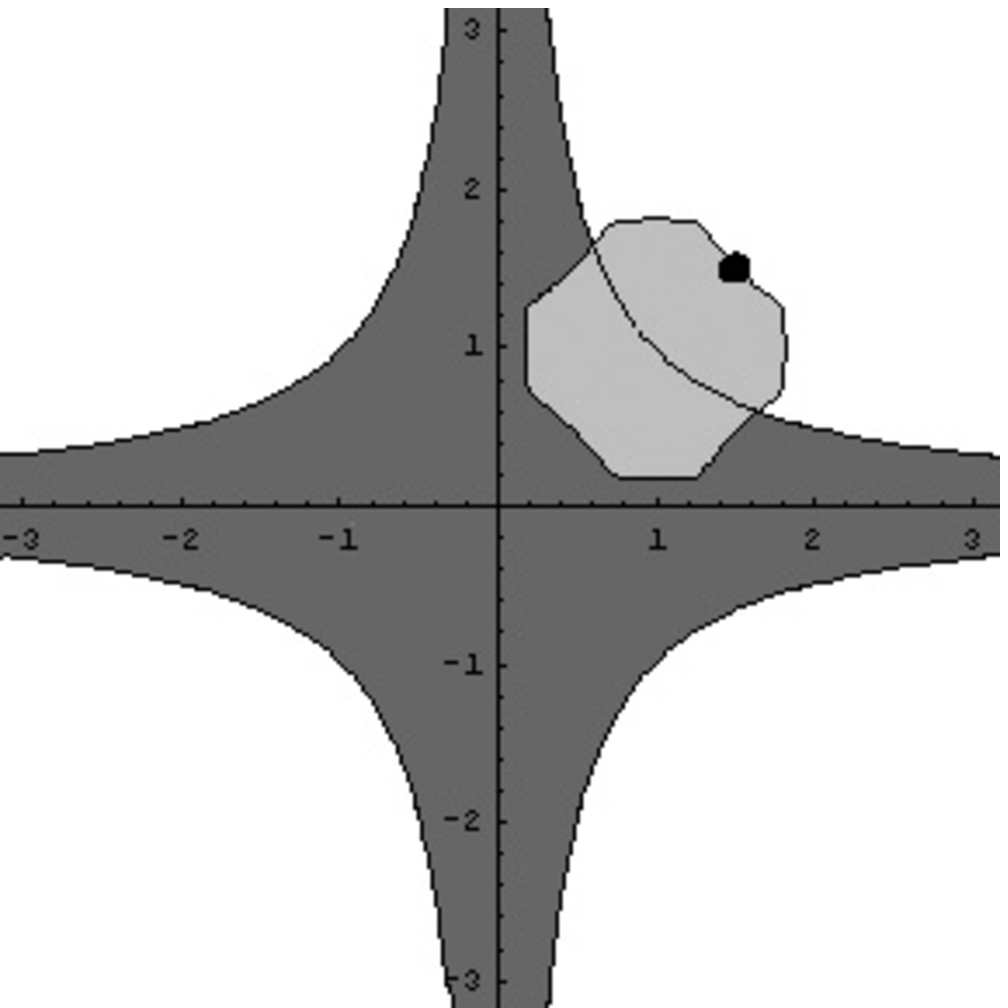} \caption{The estimate of $DA(0)$ after
2 steps for system (\ref{ex5})}
\end{figure}

\newpage

\subsubsection*{Example 6}

The following system of three difference equations is considered:

\begin{equation}
\label{ex6}
\begin{array}{lll}
\left\{\begin{array}{l}
x_{n+1}=\frac{1}{2}x_{n}y_{n}
+\frac{1}{4}x_{n}z_{n} +\frac{1}{3}x_{n}^2y_{n}+
\frac{1}{12}x_{n}^2z_{n}
     -\frac{1}{3}x_{n}y_{n}^2 -\frac{1}{12} x_{n}z_{n}^2 -\frac{1}{12} x_{n}y_{n}z_{n}\\
y_{n+1}=-\frac{1}{2}x_{n}y_{n}+\frac{1}{2}y_{n}z_{n}-\frac{1}{3}x_{n}^2y_{n}+\frac{1}{3}x_{n}y_{n}^2
  + \frac{1}{3}y_{n}^2z_{n}-\frac{1}{3} y_{n}z_{n}^2+\frac{1}{6}x_{n}y_{n}z_{n}\\
z_{n+1}=-\frac{1}{2}y_{n}z_{n}-\frac{1}{4}x_{n}z_{n}-\frac{1}{12}x_{n}^2z_{n}-\frac{1}{3}y_{n}^2z_{n}
  +\frac{1}{12}x_{n}z_{n}^2+ \frac{1}{3}y_{n}z_{n}^2-\frac{1}{12}x_{n}y_{n}z_{n}
\end{array}\right.
\end{array}
\end{equation}

The estimate of the $DA$ of the asymptotically stable fixed point
$(0,0,0)$ obtained after one step is presented in Fig. 4.

\begin{figure}[htbp]
\centering
\includegraphics*[bb=0cm 0cm 13cm
13cm, width=6cm]{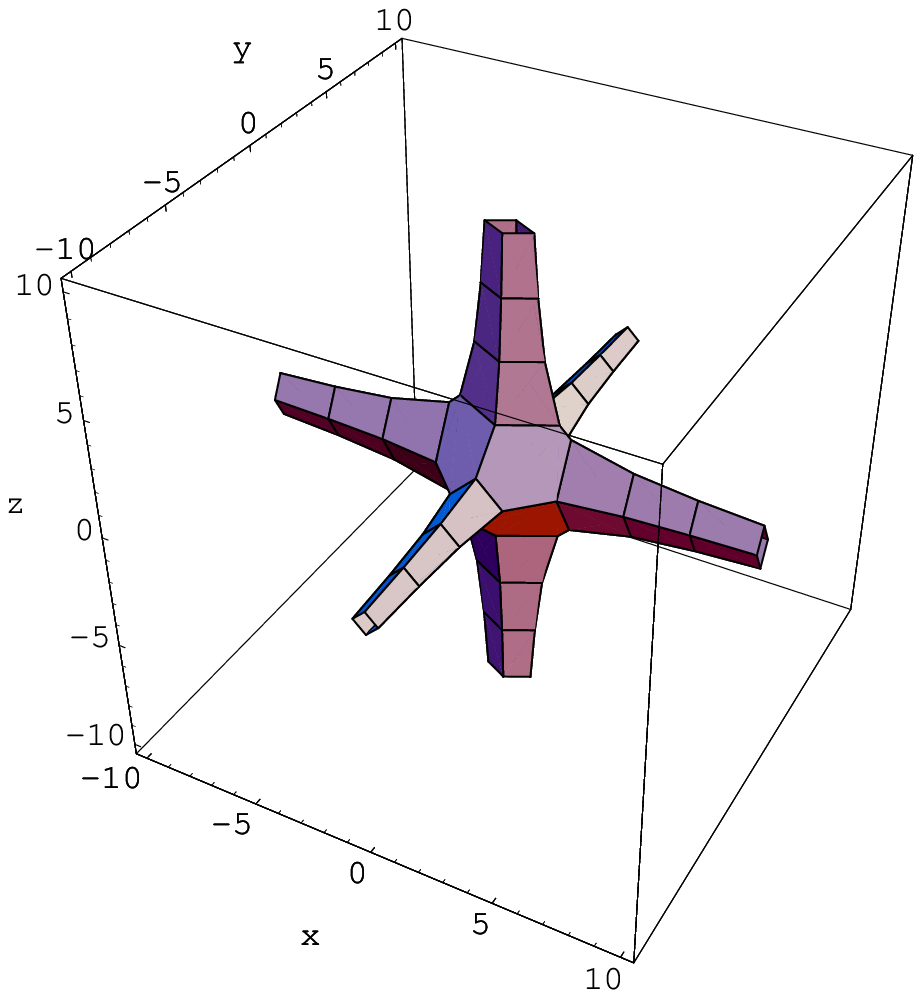} \caption{The estimate of $DA(0)$ after
1 step for system (\ref{ex6})}
\end{figure}

\newpage

\begin{center}
Table 1. Numerical data
\end{center}

\begin{tabular}{|c|c|c|c|}
  \hline
  example & order of approximation $p$ & timing for $1^{st}$ step & timing for $2^{nd}$ step \\
  \hline
  \hline
  1 & 4096 &  9.3 h & 32.5 h \\
  \hline
  2 & 164 & 10 min & 24.1 h \\
  \hline
  3 & 500 & 12.4 min & - \\
  \hline
  4 & 625 &  1.2 h & 12.4 h \\
  \hline
  5 & 256 & 16.8 min & 45.2 h \\
  \hline
  6 & 54 & 4.6 h &  - \\
  \hline
\end{tabular}

\end{document}